\documentclass{amsart}

\usepackage{amssymb}
\usepackage[all]{xy}
\usepackage{hyperref}

\newtheorem{theorem}{Theorem}[section]
\newtheorem{lemma}[theorem]{Lemma}
\newtheorem{corollary}[theorem]{Corollary}
\newtheorem{remark}[theorem]{Remark}
\theoremstyle{definition}
\newtheorem{example}[theorem]{Example}

\newcommand{\tA}{\underset{A}{\otimes}}
\newcommand{\tAe}{\underset{A^e}{\otimes}}
\newcommand{\tZ}{\underset{Z}{\otimes}}

\begin{document}

\title[Cap product]{On the cap product in Hochschild theory}

\author{Marco Armenta}

\address{CIMAT A. C., Guanajuato, M\'exico.}
\address{IMAG, Univ Montpellier, CNRS, Montpellier, France.}
\email{marco.armenta@usherbrooke.ca}


\begin{abstract}
In this paper we give an axiomatic characterization of the cap product in the Hochschild theory of associative unital algebras which are projective over a commutative unital ring. We also give an interpretation of the cap product with coefficients in the algebra via chain maps. We illustrate these results by computing the cap product for truncated polynomial algebras $k[x]/(x^N)$ and for polynomial algebras, where it is identified with the contraction of differential forms by polyvector fields.
\end{abstract}

\maketitle

\section{Introduction}

The work done in \cite{AK1} shows that the cap product is a derived invariant, as it is known for the cup product and the Gerstenhaber bracket. Moreover, the Gerstenhaber module structure of Hochschild homology over Hochschild cohomology via the combination of the cap product and the Connes differential, also known as a Tamarkin-Tsygan calculus structure, is also a derived invariant, see \cite{AK2}. Therefore, more computational techniques for the cap product are needed in order to calculate these derived invariants. We address this with two approaches. The first by following the ideas of Sanada in \cite{Sa1} and \cite{Sa2}. In his work, Sanada gives an axiomatic characterization of the cup product in Hochschild cohomology. We give an axiomatic characterization of the cap product \`a la Sanada. For the second approach we recall that the cup product can be realized by chain maps induced from cocycles, see \cite{Wi}, Section 2.1. We give an interpretation of the cap product, with coefficients in the algebra, via chain maps induced by a cocycle. We illustrate both results in Section \ref{sec:examples} by computing the cap product for the truncated polynomial algebras $k[x]/(x^N)$ and for polynomial algebras, where it is identified with the contraction of differential forms by polyvector fields.

Let $k$ be a commutative unital ring and $A$ an associative unital $k$-algebra that is projective as a $k$-module. The tensor products over $k$ will be written without subscript. We denote by $A^e$ the \emph{enveloping algebra} of $A$ which is equal to $A\otimes A^{op}$ as a $k$-module with product
\[
(a\otimes b)(a'\otimes b')=aa'\otimes b'b,
\]
for $a,a'\in A$ and $b,b'\in A^{op}$. We identify $A$-bimodules with left (and right) $A^e$-modules. The \emph{bar resolution} of $A$ will be denoted as $\mathrm{Bar}(A)_\bullet$. It is given by $\mathrm{Bar}(A)_n=A^{\otimes(n+2)}$, the differential is
\[
d_n(a_0\otimes\cdots\otimes a_{n+1})=\sum_{i=0}^{n}(-1)^ia_0\otimes\cdots\otimes a_ia_{i+1}\otimes\cdots\otimes a_{n+1},
\]
and the augmentation is the multiplication $\mu:A\otimes A\to A$ of the algebra $A$. Let $N$ and $M$ be $A^e$-modules. The \emph{Hochschild homology} of $A$ with coefficients in $N$ is defined as
\[
H_n(A,N)=H_n\big(N\tAe\mathrm{Bar}(A)_\bullet\big)
\]
for every $n\geq 0$. The \emph{Hochschild cohomology} of $A$ with coefficients in $M$ is defined as
\[
H^m(A,M)=H^m\big(Hom_{A^e}(\mathrm{Bar}(A)_\bullet,M)\big)
\]
for every $m\geq 0$. Since $A$ is $k$-projective, the bar resolution is a projective resolution of $A$ as an $A^e$-module and we have that
\[
H_n(A,N)\cong Tor_n^{A^e}(N,A)\quad\text{and}\quad H^m(A,M)\cong Ext^m_{A^e}(A,M),
\]
so that Hochschild (co)homology can be computed via any $A^e$-projective resolution of $A$. The (co)homology class of a (co)cycle will be denoted by brackets $[-]$. We write $HH_n(A)=H_n(A,A)$ and $HH^m(A)=H^m(A,A)$. The \emph{cap product} of $A$ is the graded map
\[
\cap:H_n(A,N)\otimes H^m(A,M)\to H_{n-m}(A,N\tA M)
\]
given by
\[
[x\tAe a_0\otimes\cdots\otimes a_{n+1}]\cap[t]=[x\tA t(a_0\otimes\cdots\otimes a_m\otimes 1)\tAe 1\otimes a_{m+1}\otimes\cdots\otimes a_{n+1}],
\]
for a cycle $x\tAe a_0\otimes\cdots\otimes a_{n+1}\in N\tAe A^{\otimes(n+2)}$ and a cocycle $t\in Hom_{A^e}(A^{\otimes(m+2)},M)$, where we have denoted by $1\in A$ the unit of $A$. In fact, the cap product is defined at the level of (co)chains by the same formula.

\section{Diagonal map}

Let $P_\bullet\to A$ be a projective resolution of $A$ as an $A^e$-module and denote by $d_i:P_i\to P_{i-1}$ the differential of $P_\bullet$ and the augmentation by $d_0:P_0\to A$. A \emph{diagonal map} on $P_\bullet$ is $\triangle=\{\triangle_{i,j}\}_{i,j\geq 0}$ where
\[
\triangle_{i,j}:P_{i+j}\to P_i\tA P_j
\]
is a map of $A^e$-modules such that the following equations hold:
\begin{itemize}
\item $\triangle_{i,j}d_{i+j+1}=(d_{i+1}\tA 1)\,\triangle_{i+1,j}+(-1)^i(1\tA d_{j+1})\,\triangle_{i,j+1}$,
\item $(d_0\tA d_0)\triangle_{0,0}=d_0$.
\end{itemize}
The map $\triangle_{i,j}:A^{\otimes(i+j+2)}\to A^{\otimes(i+2)}\tA A^{\otimes(j+2)}$ given by
\[
\triangle_{i,j}(a_0\otimes\cdots\otimes a_{i+j+1})=(a_0\otimes\cdots\otimes a_i\otimes 1)\tA(1\otimes a_{i+1}\otimes\cdots\otimes a_{i+j+1})
\]
defines a diagonal map on the bar resolution, see \cite{Sa2}.

\begin{remark}\label{rmk:barcap}
Let $\mathrm{Bar}(A)_\bullet\to A$ be the bar resolution and $\triangle$ the diagonal map mentioned above. Then at the level of (co)chains we have that
\[
(x\tAe a_0\otimes\cdots\otimes a_{n+1})\cap t=(\mathrm{id}\tA t\tAe\mathrm{id})\big(x\tA\triangle_{m,n-m}(a_0\otimes\cdots\otimes a_{n+1})\big),
\]
for all $x\in N$, all $a_0\otimes\cdots\otimes a_{n+1}\in A^{\otimes(n+2)}$ and all $t\in Hom_{A^e}(A^{\otimes(m+2)},M)$. Note that we have composed $\triangle$ with the quotient map
\[
P_m\tA P_{n-m}\to P_m\tAe P_{n-m}.
\]
\end{remark}

\begin{remark}\label{rmk:identify}
Let $N$, $M$ and $L$ be $A^e$-modules. We will identify $N\tA M\tAe L$ with $N\tAe M\tA L$ via the canonical isomorphism $\alpha\tA\beta\tAe\gamma\mapsto\alpha\tAe\beta\tA\gamma$.
\end{remark}

Let $\triangle:P_\bullet\to P_\bullet\otimes_AP_\bullet$ be a diagonal map. When $i=m$ and $j=n-m-1$ we have that
\begin{equation}\label{eq:diag1}
\triangle_{m,n-m-1}d_n=(d_{m+1}\tA\mathrm{id})\,\triangle_{m+1,n-m-1}+(-1)^m(\mathrm{id}\tA d_{n-m})\,\triangle_{m,n-m}\,.
\end{equation}
Each term in this equation is a map from $P_n$ to $P_m\tA P_{n-m-1}$. The term $\triangle_{m,n-m-1}d_n$ is the composition
\[
\xymatrix{P_n\ar[r]^-{d_n} & P_{n-1}\ar[rr]^-{\triangle_{m,n-m-1}} & & P_m\tA P_{n-m-1}\,,}
\]
which induces a commutative diagram
\[
\xymatrix@C=3pc{
N\tA P_n\ar[r]^-{\mathrm{id}\tA d_n}\ar[d] & N\tA P_{n-1}\ar[rr]^-{\mathrm{id}\tA\triangle_{m,n-m-1}}\ar[d] & & N\tA P_m\tA P_{n-m-1}\ar[d]\\
N\tAe P_n\ar[r]^-{\mathrm{id}\tAe d_n} & N\tAe P_{n-1}\ar[rr]^-{\mathrm{id}\tAe\triangle_{m,n-m-1}} & & N\tAe P_m\tA P_{n-m-1},
}
\]
where the vertical maps are the obvious quotient maps. The same construction for the other terms in \eqref{eq:diag1} and reordering of the factors gives
\begin{eqnarray}\label{eq:diag2}
&&(\mathrm{id}\tA\mathrm{id}\tAe d_{n-m})(\mathrm{id}\tA\triangle_{m,n-m})\\
&=&(-1)^m\,\mathrm{id}\tA(\triangle_{m,n-m-1}d_n)+(-1)^{m+1}(\mathrm{id}\tA d_{m+1}\tAe\mathrm{id})(\mathrm{id}\tA\triangle_{m+1,n-m-1}).\nonumber
\end{eqnarray}
For the sake of simplicity we will not write the quotient maps. Again, we have composed $\triangle$ with the quotient map $P_i\tA P_j\to P_i\tAe P_j$.

\section{Axiomatic characterization of the cap product}

We denote by $[N,A]$ the submodule of $N$ generated by elements of the form $ax-xa$ for $x\in N$ and $a\in A$, and by $M^A$ the submodule of $M$ formed by elements $m$ such that $am=ma$ for all $a\in A$. Let $P_\bullet\to A$ be an $A^e$-projective resolution of $A$. There is an isomorphism
\[
\begin{array}{rccc}
\psi: & H_0(A,N) & \to & N/[N,A]\\[2pt]
 & [x\tAe p] & \mapsto & [xd_0(p)],
\end{array}
\]
and recall from \cite{Sa2} that there is also an isomorphism
\[
\begin{array}{rccc}
\varphi: & H^0(A,M) & \to & M^A\\[2pt]
 & [t] & \mapsto & \widetilde{t}(1),
\end{array}
\]
where $t=\widetilde{t}d_0$ and $1\in A$ is the unit. The well-definedness of $\psi$ and $\varphi$ depends on the exactness of the sequence $P_1\to P_0\to A\to 0$. Let $Z=Z(A)$ be the center of the algebra $A$. A \emph{cap product} on $A$ is a map
\[
\cap:H_n(A,N)\otimes H^m(A,M)\to H_{n-m}(A,N\tA M)
\]
for all integers $0\leq m\leq n$ and all $A^e$-modules $N$ and $M$, that satisfies the following four axioms:
\begin{itemize}
\item[$(QI)$] $\cap$ induces a morphism of $Z$-modules
\[
\cap:H_n(A,N)\tZ H^m(A,M)\to H_{n-m}(A,N\tA M).
\]
\item[$(QII_1)$] Let
\[
\xymatrix{0\ar[r] & N_1\ar[r]^-{f} & N_2\ar[r]^-{g} & N_3\ar[r] & 0}
\]
and
\[
\xymatrix{0\ar[r] & N_1\tA M\ar[r]^-{f\otimes\mathrm{id}} & N_2\tA M\ar[r]^-{g\otimes\mathrm{id}} & N_3\tA M\ar[r] & 0}
\]
be short exact sequences of $A^e$-modules. Then
\[
\delta(\gamma\cap\epsilon)=(-1)^m(\delta\gamma)\cap\epsilon
\]
for $\gamma\in H_n(A,N_3)$ and $\epsilon\in H^m(A,M)$, where $\delta$ is the connecting homomorphism of the Hochschild homology functor $H_\bullet(A,-)$.
\item[$(QII_2)$] Let
\[
\xymatrix{0\ar[r] & M_1\ar[r]^-{f} & M_2\ar[r]^-{g} & M_3\ar[r] & 0}
\]
and
\[
\xymatrix{0\ar[r] & N\tA M_1\ar[r]^-{\mathrm{id}\otimes f} & N\tA M_2\ar[r]^-{\mathrm{id}\otimes g} & N\tA M_3\ar[r] & 0}
\]
be short exact sequences of $A^e$-modules. Then
\[
\delta(\gamma\cap\epsilon)=(-1)^{m+1}\gamma\cap(\partial\epsilon)
\]
for $\gamma\in H_n(A,N)$ and $\epsilon\in H^m(A,M_3)$, where $\partial$ is the connecting homomorphism of the Hochschild cohomology functor $H^\bullet(A,-)$.
\item[$(QIII)$] The following diagram is commutative
\[
\xymatrix{
H_0(A,N)\tZ H^0(A,M)\ar[r]^-{\cap}\ar[d]_-{\cong}^-{\psi\tZ\varphi} & H_0(A,N\tA M)\ar[d]_-{\cong}^-{\psi}\\
\dfrac{N}{[N,A]}\tZ M^A\ar[r] & \dfrac{N\tA M}{[N\tA M,A]},
}
\]
where the map in the bottom is given by $[x]\tZ y\mapsto[x\tA y]$, and the vertical maps are the isomorphisms defined at the beginning of this section.
\end{itemize}

\begin{theorem}\label{thm:cap}
There exists one and only one cap product satisfying properties $(QI)-(QIII)$.
\end{theorem}

Before proving Theorem \ref{thm:cap} we record, following \cite{Sa1} page 73, the two short exact sequences of coefficients along which we will shift dimensions, one induced and one coinduced. Let $L$ and $M$ be $A^e$-modules. We endow the $k$-modules $A\otimes L$ and $Hom_k(A,M)$ with the $A^e$-module structures
\[
a\cdot(b\otimes x)\cdot c=ab\otimes xc\qquad\text{and}\qquad(a\cdot f\cdot c)(y)=a\,f(cy),
\]
for all $a,b,c,y\in A$, $x\in L$ and $f\in Hom_k(A,M)$. There are morphisms of $A^e$-modules
\[
\pi_L:A\otimes L\to L,\quad\pi_L(b\otimes x)=bx,\qquad\iota_M:M\to Hom_k(A,M),\quad\iota_M(m)(y)=my,
\]
where $\pi_L$ is surjective and $\iota_M$ is injective. We write $K(L)=ker(\pi_L)$ and $C(M)=coker(\iota_M)$. Since $A$ is projective over $k$, each $\mathrm{Bar}(A)_n\cong A^e\otimes A^{\otimes n}$ is projective both as an $A^e$-module and as a left $A$-module, and every short exact sequence of $A^e$-modules induces short exact sequences of Hochschild chain and cochain complexes; in particular, it induces long exact sequences in $H_\bullet(A,-)$ and $H^\bullet(A,-)$.

\begin{lemma}\label{lemma:acyclic}
Let $L$ and $M$ be $A^e$-modules.
\begin{enumerate}
\item The morphism $\pi_L$ splits as a morphism of right $A$-modules via $x\mapsto 1\otimes x$, and
\[
H_j(A,A\otimes L)=0\quad\text{for all }j\geq 1.
\]
\item The morphism $\iota_M$ splits as a morphism of left $A$-modules via $f\mapsto f(1)$, and
\[
H^j(A,Hom_k(A,M))=0\quad\text{for all }j\geq 1.
\]
\end{enumerate}
\end{lemma}

\begin{proof}
(1) The map $\sigma:L\to A\otimes L$, $\sigma(x)=1\otimes x$, satisfies $\sigma(xa)=1\otimes xa=(1\otimes x)a$ and $\pi_L\sigma=\mathrm{id}_L$, so it is a right $A$-linear section of $\pi_L$.

For the vanishing, we claim that for every $A^e$-module $Y$ there is an isomorphism of $k$-modules, natural in $Y$,
\[
\Phi_Y:(A\otimes L)\tAe Y\longrightarrow L\tA Y,\qquad\Phi_Y\big((b\otimes x)\tAe y\big)=x\tA yb,
\]
where on the right hand side $L$ is regarded as a right $A$-module and $Y$ as a left $A$-module. Recall that the balancing relations in $X\tAe Y$ read $(c\,\xi\,a)\tAe y=\xi\tAe(a\,y\,c)$ for $\xi\in X$ and $a,c\in A$. For $\xi=b\otimes x$ this relation becomes $(cb\otimes xa)\tAe y=(b\otimes x)\tAe(ayc)$, and
\[
\Phi_Y\big((cb\otimes xa)\tAe y\big)=xa\tA y(cb)=x\tA aycb=\Phi_Y\big((b\otimes x)\tAe(ayc)\big),
\]
so $\Phi_Y$ is well defined. The assignment $x\tA y\mapsto(1\otimes x)\tAe y$ is balanced over $A$, since taking $c=b=1$ above gives $(1\otimes xa)\tAe y=(1\otimes x)\tAe(ay)$; it defines a map $\Psi_Y$ with $\Phi_Y\Psi_Y=\mathrm{id}$, and taking $a=1$, $c=b$ in the balancing relation gives $(b\otimes x)\tAe y=(1\otimes x)\tAe(yb)$, whence $\Psi_Y\Phi_Y=\mathrm{id}$. Naturality is immediate: for an $A^e$-linear map $g:Y\to Y'$ one has $\Phi_{Y'}\big((b\otimes x)\tAe g(y)\big)=x\tA g(y)b=x\tA g(yb)$.

Applying $\Phi$ to $Y=\mathrm{Bar}(A)_\bullet$ yields an isomorphism of complexes
\[
(A\otimes L)\tAe\mathrm{Bar}(A)_\bullet\cong L\tA\mathrm{Bar}(A)_\bullet\,.
\]
Since $A$ is $k$-projective, $\mathrm{Bar}(A)_\bullet\to A$ is also a resolution of $A$ by projective left $A$-modules: each $\mathrm{Bar}(A)_n=A\otimes A^{\otimes(n+1)}$ is a direct summand of a free left $A$-module, and the augmented complex is exact because it admits the $k$-linear contracting homotopy $a_0\otimes\cdots\otimes a_{n+1}\mapsto 1\otimes a_0\otimes\cdots\otimes a_{n+1}$. Hence
\[
H_j(A,A\otimes L)\cong H_j\big(L\tA\mathrm{Bar}(A)_\bullet\big)\cong Tor_j^A(L,A)=0\qquad\text{for all }j\geq 1,
\]
because $A$ is projective as a left $A$-module.

(2) For $f\in Hom_k(A,M)$ one has $(a\cdot f)(y)=af(y)$, so the map $r(f)=f(1)$ satisfies $r(a\cdot f)=a\,r(f)$ and $r\,\iota_M=\mathrm{id}_M$; thus $r$ is a left $A$-linear retraction of $\iota_M$.

For the vanishing, we claim that for every $A^e$-module $Y$ there is an isomorphism of $k$-modules, natural in $Y$,
\[
\Theta_Y:Hom_{A^e}\big(Y,Hom_k(A,M)\big)\longrightarrow Hom_A(Y,M),\qquad\Theta_Y(\varphi)(y)=\varphi(y)(1),
\]
where on the right hand side $Y$ and $M$ are regarded as left $A$-modules. The map $\Theta_Y(\varphi)$ is left $A$-linear because $\varphi(ay)(1)=\big(a\cdot\varphi(y)\big)(1)=a\,\varphi(y)(1)$. Conversely, given a left $A$-linear $\psi:Y\to M$, set $\varphi(y)(x)=\psi(yx)$; then $\varphi(y)\in Hom_k(A,M)$ and
\[
\varphi(ayb)(x)=\psi(aybx)=a\,\psi(ybx)=a\,\varphi(y)(bx)=\big(a\cdot\varphi(y)\cdot b\big)(x),
\]
so $\varphi$ is $A^e$-linear. The two assignments are mutually inverse: $\psi(y\cdot 1)=\psi(y)$, and in the other direction $\varphi(yx)(1)=\big(\varphi(y)\cdot x\big)(1)=\varphi(y)(x)$. Applying $\Theta$ to $Y=\mathrm{Bar}(A)_\bullet$ yields an isomorphism of cochain complexes
\[
Hom_{A^e}\big(\mathrm{Bar}(A)_\bullet,Hom_k(A,M)\big)\cong Hom_A\big(\mathrm{Bar}(A)_\bullet,M\big),
\]
and since $\mathrm{Bar}(A)_\bullet\to A$ is a projective resolution of $A$ by left $A$-modules,
\[
H^j(A,Hom_k(A,M))\cong Ext^j_A(A,M)=0\qquad\text{for all }j\geq 1,
\]
because $A$ is projective as a left $A$-module.
\end{proof}

\begin{corollary}\label{cor:shifting}
Let $N$ and $M$ be $A^e$-modules.
\begin{enumerate}
\item The sequences
\[
0\to K(N)\to A\otimes N\xrightarrow{\;\pi_N\;}N\to 0
\]
and
\[
0\to K(N)\tA M\to(A\otimes N)\tA M\xrightarrow{\;\pi_N\tA\mathrm{id}\;}N\tA M\to 0
\]
are short exact sequences of $A^e$-modules, and for every $j\geq 1$ the connecting homomorphism of the second sequence
\[
\delta:H_j(A,N\tA M)\longrightarrow H_{j-1}(A,K(N)\tA M)
\]
is injective.
\item The sequences
\[
0\to M\xrightarrow{\;\iota_M\;}Hom_k(A,M)\to C(M)\to 0
\]
and
\[
0\to N\tA M\to N\tA Hom_k(A,M)\to N\tA C(M)\to 0
\]
are short exact sequences of $A^e$-modules, and for every $j\geq 0$ the connecting homomorphism of the first sequence
\[
\partial:H^j(A,C(M))\longrightarrow H^{j+1}(A,M)
\]
is surjective.
\end{enumerate}
\end{corollary}

\begin{proof}
(1) The first sequence is exact by definition of $K(N)$. By Lemma \ref{lemma:acyclic}(1) it is split in the category of right $A$-modules, so applying the additive functor $-\tA M$ produces a split exact sequence of $k$-modules; its maps are obtained by tensoring morphisms of $A^e$-modules with $\mathrm{id}_M$, hence the second sequence is exact in the category of $A^e$-modules. The assignment $(b\otimes x)\tA m\mapsto b\otimes(x\tA m)$ defines an isomorphism of $A^e$-modules
\[
(A\otimes N)\tA M\cong A\otimes\big(N\tA M\big),
\]
so $H_j\big(A,(A\otimes N)\tA M\big)=0$ for all $j\geq 1$ by Lemma \ref{lemma:acyclic}(1). For $j\geq 1$ the long exact homology sequence of the second short exact sequence reads
\[
0=H_j\big(A,(A\otimes N)\tA M\big)\to H_j(A,N\tA M)\xrightarrow{\;\delta\;}H_{j-1}\big(A,K(N)\tA M\big),
\]
whence $\delta$ is injective.

(2) The first sequence is exact by definition of $C(M)$. By Lemma \ref{lemma:acyclic}(2) it is split in the category of left $A$-modules, so applying $N\tA-$ produces a split exact sequence of $k$-modules whose maps are $A^e$-linear; hence the second sequence is exact in the category of $A^e$-modules. For $j\geq 0$ the long exact cohomology sequence of the first sequence reads
\[
H^j(A,C(M))\xrightarrow{\;\partial\;}H^{j+1}(A,M)\to H^{j+1}\big(A,Hom_k(A,M)\big)=0,
\]
whence $\partial$ is surjective.
\end{proof}

\begin{remark}\label{rmk:K}
Under the isomorphism $(A\otimes N)\tA M\cong A\otimes(N\tA M)$ of the proof, the morphism $\pi_N\tA\mathrm{id}$ corresponds to $\pi_{N\tA M}$. Consequently $K(N)\tA M\cong K(N\tA M)$, and the second sequence in Corollary \ref{cor:shifting}(1) is canonically isomorphic to
\[
0\to K\big(N\tA M\big)\to A\otimes\big(N\tA M\big)\to N\tA M\to 0\,.
\]
\end{remark}

\begin{proof}[Proof of Theorem \ref{thm:cap}]
For a projective resolution $P_\bullet\to A$ of $A$ as an $A^e$-module and a diagonal map $\triangle:P_\bullet\to P_\bullet\tA P_\bullet$ we define
\[
\cap:N\tAe P_n\otimes Hom_{A^e}(P_m,M)\to N\tA M\tAe P_{n-m}
\]
by
\[
x\tAe p\cap t:=(\mathrm{id}\tA t\tAe\mathrm{id})\big(x\tA\triangle_{m,n-m}(p)\big).
\]
We first prove that this product descends to (co)homology. Let $x\tAe p\in N\tAe P_n$ be a cycle and let $t\in Hom_{A^e}(P_m,M)$ be a cocycle. By \eqref{eq:diag2} we have that
\begin{eqnarray*}
&&(\mathrm{id}\tA\mathrm{id}\tAe d_{n-m})(x\tAe p\cap t)\\
&=&(\mathrm{id}\tA\mathrm{id}\tAe d_{n-m})(\mathrm{id}\tA t\tAe\mathrm{id})(x\tA\triangle_{m,n-m}(p))\\
&=&(\mathrm{id}\tA t\tAe\mathrm{id})(\mathrm{id}\tA\mathrm{id}\tAe d_{n-m})(x\tA\triangle_{m,n-m}(p))\\
&=&(-1)^m(\mathrm{id}\tA t\tAe\mathrm{id})(x\tA\triangle_{m,n-m-1}\big(d_n(p)\big))\\
&&\quad+\,(-1)^{m+1}(\mathrm{id}\tA t\tAe\mathrm{id})(\mathrm{id}\tA d_{m+1}\tAe\mathrm{id})(x\tA\triangle_{m+1,n-m-1}(p)).
\end{eqnarray*}
By Remark \ref{rmk:identify} and the discussion following it, we have that
\[
(\mathrm{id}\tA\triangle_{m,n-m-1})(\mathrm{id}\tA d_n)(x\tAe p)=(\mathrm{id}\tAe\triangle_{m,n-m-1})(\mathrm{id}\tAe d_n)(x\tAe p)=0.
\]
Observe also that $(\mathrm{id}\tA t\tAe\mathrm{id})(\mathrm{id}\tA d_{m+1}\tAe\mathrm{id})(x\tA\triangle_{m+1,n-m-1}(p))$ is equal to
\[
(\mathrm{id}\tA td_{m+1}\tAe\mathrm{id})(x\tA\triangle_{m+1,n-m-1}(p))=0
\]
since $t$ is a cocycle. Therefore $x\tAe p\cap t$ is a cycle. We now compose equation \eqref{eq:diag2} on the left with $(\mathrm{id}\tA t\tAe\mathrm{id})$ to obtain
\[
(\mathrm{id}\tA\mathrm{id}\tAe d_{n-m})(x\tAe p\cap t)=(-1)^mx\tAe d_n(p)\cap t+(-1)^{m+1}x\tAe p\cap(td_{m+1}).
\]
From this equation it is clear that $\cap$ is well-defined in (co)homology
\[
\cap:H_n(A,N)\otimes H^m(A,M)\to H_{n-m}(A,N\tA M).
\]
We now prove that this product satisfies the $Q$-properties. It is clear that it satisfies property $(QI)$. Take short exact sequences as in $(QII_1)$. Let $\gamma$ be an element of $H_n(A,N_3)$ represented by a cycle $x\tAe p$ for $x\in N_3$ and $p\in P_n$. Let $\epsilon$ be an element of $H^m(A,M)$ represented by a cocycle $t\in Hom_{A^e}(P_m,M)$. We use the snake lemma to construct $\delta(\gamma)\in H_{n-1}(A,N_1)$. We have a diagram
\[
\xymatrix{
 & N_2\tAe P_n\ar[r]^-{g\tAe\mathrm{id}}\ar[d]^-{\mathrm{id}\tAe d_n} & N_3\tAe P_n\ar[r] & 0\\
N_1\tAe P_{n-1}\ar[r]^-{f\tAe\mathrm{id}} & N_2\tAe P_{n-1}. & &
}
\]
There exists $y=\sum_iy_i\tAe p_i\in N_2\tAe P_n$ such that $x\tAe p=\sum_ig(y_i)\tAe p_i$. There also exists $z=\sum_jz_j\tAe q_j\in N_1\tAe P_{n-1}$ such that $\sum_jf(z_j)\tAe q_j=\sum_iy_i\tAe d_n(p_i)$. Therefore $\delta(\gamma)=[z]\in H_{n-1}(A,N_1)$. We now compute $\delta(\gamma\cap\epsilon)$. There is a diagram
\[
\xymatrix{
 & N_2\tA M\tAe P_{n-m}\ar[r]^-{g\tA\mathrm{id}\tAe\mathrm{id}}\ar[d]^-{\mathrm{id}\tA\mathrm{id}\tAe d_{n-m}} & N_3\tA M\tAe P_{n-m}\ar[r] & 0\\
N_1\tA M\tAe P_{n-m-1}\ar[r]^-{f\tA\mathrm{id}\tAe\mathrm{id}} & N_2\tA M\tAe P_{n-m-1}. & &
}
\]
The element $x\tAe p\cap t$ lies in $N_3\tA M\tAe P_{n-m}$ and $y\cap t$ is an element of $N_2\tA M\tAe P_{n-m}$. The morphism $g\tA\mathrm{id}\tAe\mathrm{id}$ maps $y\cap t$ to $x\tAe p\cap t$. Indeed,
\begin{eqnarray*}
x\tAe p\cap t&=&\Big(\sum_ig(y_i)\tAe p_i\Big)\cap t\\
&=&\sum_i\Big(g(y_i)\tAe p_i\cap t\Big)\\
&=&\sum_i(\mathrm{id}\tA t\tAe\mathrm{id})(g(y_i)\tA\triangle_{m,n-m}(p_i))\\
&=&\sum_i(\mathrm{id}\tA t\tAe\mathrm{id})(g\tA\mathrm{id}\tAe\mathrm{id})(y_i\tA\triangle_{m,n-m}(p_i))\\
&=&(g\tA\mathrm{id}\tAe\mathrm{id})\sum_i(\mathrm{id}\tA t\tAe\mathrm{id})(y_i\tA\triangle_{m,n-m}(p_i))\\
&=&(g\tA\mathrm{id}\tAe\mathrm{id})\Big(\sum_i(y_i\tAe p_i)\cap t\Big)\\
&=&(g\tA\mathrm{id}\tAe\mathrm{id})(y\cap t).
\end{eqnarray*}
Since $t$ is a cocycle, by \eqref{eq:diag2} we get that
\[
(\mathrm{id}\tA t\tAe\mathrm{id})(\mathrm{id}\tA\triangle_{m,n-m-1}\,d_n)=(-1)^m(\mathrm{id}\tA t\tAe d_{n-m})(\mathrm{id}\tA\triangle_{m,n-m}).
\]
Observe also that $z\cap t$ is an element of $N_1\tA M\tAe P_{n-m-1}$ and then
\begin{eqnarray*}
(f\tA\mathrm{id}\tAe\mathrm{id})(z\cap t)&=&(f\tA\mathrm{id}\tAe\mathrm{id})\Big(\sum_jz_j\tAe q_j\cap t\Big)\\
&=&(f\tA\mathrm{id}\tAe\mathrm{id})\sum_j(\mathrm{id}\tA t\tAe\mathrm{id})(z_j\tA\triangle_{m,n-m-1}(q_j))\\
&=&\sum_j(\mathrm{id}\tA t\tAe\mathrm{id})(f\tA\mathrm{id}\tAe\mathrm{id})(z_j\tA\triangle_{m,n-m-1}(q_j))\\
&=&\sum_j(\mathrm{id}\tA t\tAe\mathrm{id})(f(z_j)\tA\triangle_{m,n-m-1}(q_j))\\
&=&(\mathrm{id}\tA t\tAe\mathrm{id})\Big(\sum_jf(z_j)\tA\triangle_{m,n-m-1}(q_j)\Big)\\
&=&\Big(\sum_jf(z_j)\tAe q_j\Big)\cap t\\
&=&\Big(\sum_iy_i\tAe d_n(p_i)\Big)\cap t\\
&=&\sum_i\Big(y_i\tAe d_n(p_i)\cap t\Big)\\
&=&\sum_i(\mathrm{id}\tA t\tAe\mathrm{id})(y_i\tA\triangle_{m,n-m-1}\big(d_n(p_i)\big))\\
&=&(-1)^m\sum_i(\mathrm{id}\tA t\tAe d_{n-m})(y_i\tA\triangle_{m,n-m}(p_i))\\
&=&(-1)^m(\mathrm{id}\tA\mathrm{id}\tAe d_{n-m})\sum_i(\mathrm{id}\tA t\tAe\mathrm{id})(y_i\tA\triangle_{m,n-m}(p_i))\\
&=&(-1)^m(\mathrm{id}\tA\mathrm{id}\tAe d_{n-m})\Big(\sum_iy_i\tAe p_i\cap t\Big)\\
&=&(-1)^m(\mathrm{id}\tA\mathrm{id}\tAe d_{n-m})(y\cap t).
\end{eqnarray*}
Therefore
\begin{eqnarray*}
\delta(\gamma\cap\epsilon)&=&(-1)^m[z\cap t]\\
&=&(-1)^m[z]\cap[t]\\
&=&(-1)^m(\delta\gamma)\cap\epsilon.
\end{eqnarray*}
Now take short exact sequences as in $(QII_2)$. Let $\gamma\in H_n(A,N)$ be represented by a cycle $x\tAe p\in N\tAe P_n$ and let $\epsilon\in H^m(A,M_3)$ be represented by a cocycle $t_3\in Hom_{A^e}(P_m,M_3)$. We use the snake lemma to construct $\partial(\epsilon)$. There is a diagram
\[
\xymatrix{
 & Hom_{A^e}(P_m,M_2)\ar[r]^-{g_\#}\ar[d]^-{d^\#_{m+1}} & Hom_{A^e}(P_m,M_3)\ar[r] & 0\\
Hom_{A^e}(P_{m+1},M_1)\ar[r]^-{f_\#} & Hom_{A^e}(P_{m+1},M_2). & &
}
\]
Where $f_\#$ denotes the image of $f$ under the covariant $Hom$ functor, and $d^\#_{m+1}$ the image of $d_{m+1}$ under the contravariant $Hom$ functor. There exists $t_2\in Hom_{A^e}(P_m,M_2)$ such that $g_\#t_2=gt_2=t_3$. There also exists $t_1\in Hom_{A^e}(P_{m+1},M_1)$ such that $ft_1=f_\#t_1=d^\#_{m+1}t_2=t_2d_{m+1}$. Therefore $\partial(\epsilon)=[t_1]\in H^{m+1}(A,M_1)$. We now compute $\delta(\gamma\cap\epsilon)$. There is a diagram
\[
\xymatrix{
 & N\tA M_2\tAe P_{n-m}\ar[r]^-{\mathrm{id}\tA g\tAe\mathrm{id}}\ar[d]^-{\mathrm{id}\tA\mathrm{id}\tAe d_{n-m}} & N\tA M_3\tAe P_{n-m}\ar[r] & 0\\
N\tA M_1\tAe P_{n-m-1}\ar[r]^-{\mathrm{id}\tA f\tAe\mathrm{id}} & N\tA M_2\tAe P_{n-m-1}. & &
}
\]
The element $x\tAe p\cap t_3$ lies in $N\tA M_3\tAe P_{n-m}$ and $x\tAe p\cap t_2$ lies in $N\tA M_2\tAe P_{n-m}$. Then
\begin{eqnarray*}
(\mathrm{id}\tA g\tAe\mathrm{id})(x\tAe p\cap t_2)&=&(\mathrm{id}\tA g\tAe\mathrm{id})(\mathrm{id}\tA t_2\tAe\mathrm{id})\Big(x\tA\triangle_{m,n-m}(p)\Big)\\
&=&(\mathrm{id}\tA gt_2\tAe\mathrm{id})\Big(x\tA\triangle_{m,n-m}(p)\Big)\\
&=&(\mathrm{id}\tA t_3\tAe\mathrm{id})\Big(x\tA\triangle_{m,n-m}(p)\Big)\\
&=&x\tAe p\cap t_3.
\end{eqnarray*}
Recall from the first part of the proof that
\[
(\mathrm{id}\tA d_{m+1}\tAe\mathrm{id})\big(x\tA\triangle_{m+1,n-m-1}(p)\big)=(-1)^{m+1}(\mathrm{id}\tA\mathrm{id}\tAe d_{n-m})(x\tA\triangle_{m,n-m}(p)),
\]
since $x\tAe p$ is a cycle. Now $x\tAe p\cap t_1$ is an element of $N\tA M_1\tAe P_{n-m-1}$ and we get
\begin{eqnarray*}
&&(\mathrm{id}\tA f\tAe\mathrm{id})(x\tAe p\cap t_1)\\
&=&(\mathrm{id}\tA f\tAe\mathrm{id})(\mathrm{id}\tA t_1\tAe\mathrm{id})\Big(x\tA\triangle_{m+1,n-m-1}(p)\Big)\\
&=&(\mathrm{id}\tA ft_1\tAe\mathrm{id})\Big(x\tA\triangle_{m+1,n-m-1}(p)\Big)\\
&=&(\mathrm{id}\tA t_2d_{m+1}\tAe\mathrm{id})\Big(x\tA\triangle_{m+1,n-m-1}(p)\Big)\\
&=&(\mathrm{id}\tA t_2\tAe\mathrm{id})(\mathrm{id}\tA d_{m+1}\tAe\mathrm{id})\Big(x\tA\triangle_{m+1,n-m-1}(p)\Big)\\
&=&(-1)^{m+1}(\mathrm{id}\tA t_2\tAe\mathrm{id})(\mathrm{id}\tA\mathrm{id}\tAe d_{n-m})\big(x\tA\triangle_{m,n-m}(p)\big)\\
&=&(-1)^{m+1}(\mathrm{id}\tA\mathrm{id}\tAe d_{n-m})(\mathrm{id}\tA t_2\tAe\mathrm{id})\big(x\tA\triangle_{m,n-m}(p)\big)\\
&=&(-1)^{m+1}(\mathrm{id}\tA\mathrm{id}\tAe d_{n-m})\Big(x\tAe p\cap t_2\Big).
\end{eqnarray*}
Therefore
\begin{eqnarray*}
\delta(\gamma\cap\epsilon)&=&(-1)^{m+1}[x\tAe p\cap t_1]\\
&=&(-1)^{m+1}[x\tAe p]\cap[t_1]\\
&=&(-1)^{m+1}\gamma\cap(\partial\epsilon).
\end{eqnarray*}
To prove $(QIII)$, take $[x\tAe p]\in H_0(A,N)$ and $[t]\in H^0(A,M)$. Denote $\triangle_{0,0}(p)=\triangle^1_{0,0}(p)\tA\triangle^2_{0,0}(p)$. Then
\begin{eqnarray*}
\psi([x\tAe p]\cap[t])&=&\psi[x\tA t\triangle^1_{0,0}(p)\tAe\triangle^2_{0,0}(p)]\\
&=&[x\tA\widetilde{t}d_0\triangle^1_{0,0}(p)d_0\triangle^2_{0,0}(p)]\\
&=&[x\tA\widetilde{t}(1)(d_0\tA d_0)\triangle_{0,0}(p)]\\
&=&[x\tA\widetilde{t}(1)d_0(p)],
\end{eqnarray*}
while
\begin{eqnarray*}
(\psi\tZ\varphi)\big([x\tAe p]\tZ[t]\big)&=&[xd_0(p)]\tZ\widetilde{t}(1)\\
&=&[d_0(p)x]\tZ\widetilde{t}(1)\;\mapsto\;[d_0(p)x\tA\widetilde{t}(1)]\\
&&\phantom{[d_0(p)x]\tZ\widetilde{t}(1)\;\mapsto\;\,}=[x\tA\widetilde{t}(1)d_0(p)].
\end{eqnarray*}
We now prove uniqueness. Let $\cap$ and $\cap'$ be two cap products satisfying axioms $(QI)$, $(QII_1)$, $(QII_2)$ and $(QIII)$. We say that $\cap$ and $\cap'$ \emph{agree in degree $(n,m)$} if $\alpha\cap\beta=\alpha\cap'\beta$ for all $A^e$-modules $N$ and $M$ and all $\alpha\in H_n(A,N)$, $\beta\in H^m(A,M)$. By $(QI)$ and $(QIII)$, the two products agree in degree $(0,0)$. We prove two induction steps.

\medskip
\noindent\emph{Step 1: if the products agree in degree $(n,m)$ with $m<n$, then they agree in degree $(n,m+1)$.} Let $\alpha\in H_n(A,N)$ and $\beta\in H^{m+1}(A,M)$. By Corollary \ref{cor:shifting}(2), the pair of short exact sequences
\[
0\to M\to Hom_k(A,M)\to C(M)\to 0
\]
and
\[
0\to N\tA M\to N\tA Hom_k(A,M)\to N\tA C(M)\to 0
\]
satisfies the hypotheses of $(QII_2)$, and $\partial:H^m(A,C(M))\to H^{m+1}(A,M)$ is surjective; choose $\beta'\in H^m(A,C(M))$ with $\partial\beta'=\beta$. The elements $\alpha\cap\beta'$ and $\alpha\cap'\beta'$ coincide since they are in degree $(n,m)$. Then, by $(QII_2)$,
\[
\alpha\cap\beta=\alpha\cap(\partial\beta')=(-1)^{m+1}\delta(\alpha\cap\beta')=(-1)^{m+1}\delta(\alpha\cap'\beta')=\alpha\cap'(\partial\beta')=\alpha\cap'\beta\,.
\]

\medskip
\noindent\emph{Step 2: if the products agree in degree $(n,m)$ with $m\leq n$, then they agree in degree $(n+1,m)$.} Let $\alpha\in H_{n+1}(A,N)$ and $\beta\in H^m(A,M)$. By Corollary \ref{cor:shifting}(1), the pair of short exact sequences
\[
0\to K(N)\to A\otimes N\to N\to 0
\]
and
\[
0\to K(N)\tA M\to(A\otimes N)\tA M\to N\tA M\to 0
\]
satisfies the hypotheses of $(QII_1)$. Let $\delta\alpha\in H_n(A,K(N))$ be the image of $\alpha$ under the connecting homomorphism of the first sequence. The elements $(\delta\alpha)\cap\beta$ and $(\delta\alpha)\cap'\beta$ coincide since they are in degree $(n,m)$. By $(QII_1)$,
\[
\delta(\alpha\cap\beta)=(-1)^m(\delta\alpha)\cap\beta=(-1)^m(\delta\alpha)\cap'\beta=\delta(\alpha\cap'\beta),
\]
where now $\delta:H_{n-m+1}(A,N\tA M)\to H_{n-m}\big(A,K(N)\tA M\big)$ is the connecting homomorphism of the second sequence. Since $n-m+1\geq 1$, this $\delta$ is injective by Corollary \ref{cor:shifting}(1), and therefore $\alpha\cap\beta=\alpha\cap'\beta$.

\medskip
Starting from degree $(0,0)$ and applying Step 2 repeatedly, the products agree in every degree $(n,0)$ with $n\geq 0$; for each fixed $n$, applying Step 1 repeatedly they agree in the degrees $(n,1),\dots,(n,n)$. Hence $\cap$ and $\cap'$ agree in every degree $(n,m)$ with $0\leq m\leq n$, which proves uniqueness.
\end{proof}

\section{The cap product via chain maps}

We proceed to give an interpretation of the cap product in terms of chain maps when we take $N=M=A$, that is
\[
\cap:HH_n(A)\otimes HH^m(A)\to HH_{n-m}(A),
\]
for all integers $0\leq m\leq n$. This is the version of the cap product used for the Tamarkin-Tsygan calculus structure in Hochschild (co)homology, see \cite{AK2}. Let $P_\bullet\to A$ be an $A^e$-projective resolution of $A$ and let $t\in Hom_{A^e}(P_m,A)$ be a cocycle. We proceed to construct a chain map as in \cite{Wi}, Section 2.1. Define $K_m=Im(d_m)=ker(d_{m-1})\subseteq P_{m-1}$. Since $t$ is a cocycle there is a map $\widehat{t}:K_m\to A$ such that $\widehat{t}\,\epsilon_m=t$, where $\epsilon_m:P_m\to K_m$ is the corestriction of $d_m$, which induces an isomorphism $K_m\cong P_m/Im(d_{m+1})$. By the comparison theorem, see \cite{Ro} pages 340-341, there exist maps $t_i:P_{m+i}\to P_i$ for each $i\geq 0$ such that the following is a commutative diagram with exact rows
\[
\xymatrix{
\cdots\ar[r]^-{d_{m+i+1}} & P_{m+i}\ar[r]^-{d_{m+i}}\ar[d]^-{t_i} & \cdots\ar[r]^-{d_{m+1}} & P_m\ar[r]^-{\epsilon_m}\ar[d]^-{t_0}\ar[dr]^(0.4){t} & K_m\ar[r]\ar[d]^-{\widehat{t}} & 0\\
\cdots\ar[r]^-{d_{i+1}} & P_i\ar[r]^-{d_i} & \cdots\ar[r]^-{d_1} & P_0\ar[r]^-{d_0} & A\ar[r] & 0.
}
\]
Let $a\tAe p\in A\tAe P_n$ be a cycle. We define an element of $A\tAe P_{n-m}$ by
\[
(a\tAe p)\,\widetilde{\cap}\,t:=a\tAe t_{n-m}(p).
\]
The lift $t_\bullet:P_{m+\bullet}\to P_\bullet$ is unique up to homotopy and then for any other lift $t'_\bullet:P_{m+\bullet}\to P_\bullet$ the elements $a\tAe t_{n-m}(p)$ and $a\tAe t'_{n-m}(p)$ represent the same homology class. Observe also that $a\tAe t_{n-m}(p)$ is a cycle since $a\tAe p$ is a cycle and $t_\bullet$ is a chain map. If $a\tAe p$ is a boundary, then $a\tAe t_{n-m}(p)$ is a boundary since $t_\bullet$ is a chain map. If $t$ is a coboundary, say $t=sd_m$ with $s\in Hom_{A^e}(P_{m-1},A)$, take $\widetilde{s}:P_{m-1}\to P_0$ with $d_0\widetilde{s}=s$, which exists since $P_{m-1}$ is projective and $d_0$ is surjective. Then $t_0=\widetilde{s}\,d_m$ and $t_j=0$ for all $j\geq 1$ is a lift of $\widehat{t}$, see \cite{Wi}, Section 2.1, and then $(a\tAe p)\,\widetilde{\cap}\,t=0$: this is clear for $n>m$, and for $n=m$ we have
\[
a\tAe\widetilde{s}\big(d_m(p)\big)=(\mathrm{id}\tAe\widetilde{s})\big(a\tAe d_m(p)\big)=0
\]
since $a\tAe p$ is a cycle. Therefore the $\widetilde{\cap}$-product induces a well-defined product in (co)homology.

\begin{remark}\label{rmk:resolution}
The class $\big[a\tAe t_{n-m}(p)\big]$ is also independent of the chosen resolution. Let $P'_\bullet\to A$ be another $A^e$-projective resolution of $A$ with differential $d'_\bullet$, and let $\Phi_\bullet:P'_\bullet\to P_\bullet$ and $\Psi_\bullet:P_\bullet\to P'_\bullet$ be comparison chain maps over the identity of $A$. Then $t\Phi_m\in Hom_{A^e}(P'_m,A)$ is a cocycle representing the same cohomology class as $t$, and $\Psi_\bullet\,t_{m+\bullet}\,\Phi_{m+\bullet}$ is a lift of $\widehat{t\Phi_m}$: for $j\geq 1$,
\begin{eqnarray*}
d'_j\Psi_jt_j\Phi_{m+j}&=&\Psi_{j-1}d_jt_j\Phi_{m+j}=\Psi_{j-1}t_{j-1}d_{m+j}\Phi_{m+j}\\
&=&\Psi_{j-1}t_{j-1}\Phi_{m+j-1}d'_{m+j},
\end{eqnarray*}
and $d'_0\Psi_0t_0\Phi_m=d_0t_0\Phi_m=t\Phi_m$. Since $A\tAe\Phi_\bullet$ and $A\tAe\Psi_\bullet$ are mutually inverse in homology, the class $\big[a\tAe t_{n-m}(p)\big]$ depends only on the homology class of $a\tAe p$ and on the cohomology class of $t$, and not on the resolution or on the lift.
\end{remark}

We now compare the $\widetilde{\cap}$-product with the $\cap$-product.

\begin{theorem}\label{thm:chainmaps}
Let $P_\bullet\to A$ be an $A^e$-projective resolution of $A$ and let $t\in Hom_{A^e}(P_m,A)$ be a cocycle with a chain map $t_\bullet:P_{m+\bullet}\to P_\bullet$ as above. Then, for every cycle $a\tAe p\in A\tAe P_n$ and all integers $0\leq m\leq n$,
\[
[a\tAe p]\cap[t]=(-1)^{m(n-m)}\big[a\tAe t_{n-m}(p)\big].
\]
In particular, $\widetilde{\cap}=(-1)^{m(n-m)}\cap$ on $HH_n(A)\otimes HH^m(A)$, and the two products coincide whenever $m$ or $n-m$ is even.
\end{theorem}

\begin{proof}
By Remark \ref{rmk:resolution} and the paragraph preceding it, both sides of the displayed formula depend only on the classes $[a\tAe p]$ and $[t]$, and are independent of the resolution and of the lift. It therefore suffices to prove the formula on the bar resolution, for one lift of each cocycle. Let $t\in Hom_{A^e}(A^{\otimes(m+2)},A)$ be a cocycle and define maps $s_i:A^{\otimes(m+i+2)}\to A^{\otimes(i+2)}$ by
\[
s_i(a_0\otimes\cdots\otimes a_{m+i+1})=t(a_0\otimes\cdots\otimes a_m\otimes 1)\otimes a_{m+1}\otimes\cdots\otimes a_{m+i+1}.
\]
We claim that
\begin{equation}\label{eq:signsquare}
d_i\,s_i=(-1)^m\,s_{i-1}\,d_{m+i},\qquad i\geq 1.
\end{equation}
Indeed, in $s_{i-1}d_{m+i}(a_0\otimes\cdots\otimes a_{m+i+1})$ the summands of the bar differential with $0\leq j\leq m$ leave $a_{m+2},\dots,a_{m+i+1}$ untouched and feed the contracted tuple to $t$; by the cocycle identity
\[
0=t\big(d_{m+1}(a_0\otimes\cdots\otimes a_{m+1}\otimes 1)\big)
\]
and right $A$-linearity of $t$ they assemble to
\[
(-1)^m\,t(a_0\otimes\cdots\otimes a_m\otimes 1)\,a_{m+1}\otimes a_{m+2}\otimes\cdots\otimes a_{m+i+1},
\]
which is $(-1)^m$ times the first summand of $d_is_i(a_0\otimes\cdots\otimes a_{m+i+1})$. The summands with $m+1\leq j\leq m+i$ are $(-1)^m$ times the remaining summands of $d_is_i$, under the index shift $r=j-m$. This proves \eqref{eq:signsquare}. Moreover,
\[
d_0s_0(a_0\otimes\cdots\otimes a_{m+1})=t(a_0\otimes\cdots\otimes a_m\otimes 1)\,a_{m+1}=t(a_0\otimes\cdots\otimes a_{m+1}),
\]
by right $A$-linearity of $t$. Hence $t_i:=(-1)^{mi}s_i$ is a lift of $\widehat{t}$, and for a cycle $a\tAe(a_0\otimes\cdots\otimes a_{n+1})$ we get
\begin{eqnarray*}
&&a\tAe(a_0\otimes\cdots\otimes a_{n+1})\,\widetilde{\cap}\,t\;=\;a\tAe t_{n-m}(a_0\otimes\cdots\otimes a_{n+1})\\
&=&(-1)^{m(n-m)}\,a\tAe t(a_0\otimes\cdots\otimes a_m\otimes 1)\otimes a_{m+1}\otimes\cdots\otimes a_{n+1}\\
&=&(-1)^{m(n-m)}\,at(a_0\otimes\cdots\otimes a_m\otimes 1)\tAe 1\otimes a_{m+1}\otimes\cdots\otimes a_{n+1}\\
&=&(-1)^{m(n-m)}\,\Big(a\tAe(a_0\otimes\cdots\otimes a_{n+1})\Big)\cap t.
\end{eqnarray*}
\end{proof}

\section{Examples}\label{sec:examples}

We compute the cap product of Theorem \ref{thm:cap} in two families of examples, working on small resolutions through the chain maps of Section 4. The two families sit at opposite extremes: for the smooth algebra $k[x_1,\dots,x_d]$ the cap action of $HH^\bullet(A)$ on $HH_\bullet(A)$ is the full contraction calculus of polyvector fields on differential forms, while for $k[x]/(x^N)$ it is multiplication combined with the periodicity operator of the $2$-periodic resolution, except on the quadrant of even homology and odd cohomology, where it degenerates to a single obstruction class, nonzero over a field exactly in characteristic $2$ with $N\equiv 2\pmod 4$.

Throughout this section we keep the notation of Section 4: for a cocycle $t\in Hom_{A^e}(P_m,A)$ we write $\widehat{t}:K_m\to A$ for the induced map on $K_m=Im(d_m)\subseteq P_{m-1}$, and we call a family of maps $T_\bullet:P_{m+\bullet}\to P_\bullet$ a \emph{lift of} $\widehat{t}$ if $d_0T_0=t$ and $d_jT_j=T_{j-1}d_{m+j}$ for all $j\geq 1$. By Theorem \ref{thm:chainmaps}, for every lift $T_\bullet$, every cycle $a\tAe p\in A\tAe P_n$ and all integers $0\leq m\leq n$,
\[
[a\tAe p]\cap[t]=(-1)^{m(n-m)}\big[a\tAe T_{n-m}(p)\big].
\]

\begin{example}[Truncated polynomial algebras]\label{ex:trunc}
Let $N\geq 2$, let $k$ be any commutative unital ring and $A=k[x]/(x^N)$, which is free as a $k$-module. Write $B=A^e\cong k[y,z]/(y^N,z^N)$ with $y=x\otimes 1$, $z=1\otimes x$, and
\[
u=y-z,\qquad v=\sum_{i+j=N-1}y^iz^j\;\in B,\qquad uv=y^N-z^N=0.
\]
Then
\[
\cdots\xrightarrow{\;u\;}B\xrightarrow{\;v\;}B\xrightarrow{\;u\;}B\xrightarrow{\;d_0\;}A\longrightarrow 0,
\]
with $d_{2i-1}=u\cdot$, $d_{2i}=v\cdot$ for $i\geq 1$ and $d_0$ the multiplication, is an $A^e$-free resolution of $A$. Exactness over any $k$ is a direct check on the $k$-basis $\{y^iz^j\}_{0\leq i,j\leq N-1}$: on each slice of fixed total degree, $u\cdot$ acts by consecutive differences of basis vectors, a triangular system, while $v\cdot y^iz^j=\sum_{p+q=N-1+i+j,\;p,q\leq N-1}y^pz^q$ depends only on $i+j$ and vanishes for $i+j\geq N$; telescoping then identifies $ker(u\cdot)=Im(v\cdot)$, $ker(v\cdot)=Im(u\cdot)$ and $ker(d_0)=Im(u\cdot)$.

Under $A\tAe B\cong A\cong Hom_{A^e}(B,A)$ the maps induced by $u\cdot$ and $v\cdot$ are multiplication by $0$ and by $Nx^{N-1}$, so for $i\geq 1$
\[
HH_0(A)=A,\qquad HH_{2i-1}(A)=A/Nx^{N-1}A,\qquad HH_{2i}(A)=Ann_A(Nx^{N-1}),
\]
\[
HH^0(A)=A,\qquad HH^{2i-1}(A)=Ann_A(Nx^{N-1}),\qquad HH^{2i}(A)=A/Nx^{N-1}A.
\]
If $N$ is invertible in $k$ these are $k[x]/(x^{N-1})$ and $(x)\cong k[x]/(x^{N-1})$; if $N=0$ in $k$ they are all equal to $A$. A class of $HH_n(A)$ is represented by a cycle $a\tAe 1_B$ with $a\in A$, and a class of $HH^m(A)$ by the $A^e$-linear map $t_\lambda:B\to A$, $t_\lambda(\xi)=\mu(\xi)\lambda$, with $\lambda\in A$; for odd $m$ the cocycle condition reads $Nx^{N-1}\lambda=0$. We write $[a]_n$ for the class of $a\tAe 1_B$ in $HH_n(A)$ and $[\lambda]_m$ for the class of $t_\lambda$ in $HH^m(A)$. Via the comparison $\Phi_1:\mathrm{Bar}(A)_1\to B$, $\Phi_1(1\otimes a\otimes 1)=q_a$ with $a\otimes 1-1\otimes a=q_au$, one has $\mu(q_{x^r})=rx^{r-1}$, so $[\lambda]_1\in HH^1(A)$ is the class of the derivation $D$ determined by $D(x)=\lambda$.

\emph{Lifts.} Fix $\lambda=\sum_lc_lx^l$ and the lift $\Lambda=\sum_lc_lz^l\in B$. If $m$ is even, then $T_j:=\Lambda\cdot$ for all $j\geq 0$ is a lift of $\widehat{t_\lambda}$: every square commutes because $B$ is commutative, and $d_0T_0=t_\lambda$. If $m$ is odd, the first square demands $u\,T_1(\xi)=\Lambda v\xi$, that is, an element $\Theta\in B$ with $u\Theta=v\Lambda$. Since $ker(d_0)=uB$, such a $\Theta$ exists if and only if $\mu(v\Lambda)=Nx^{N-1}\lambda$ vanishes, precisely the cocycle condition, and then
\[
T_{2j}=\Lambda\cdot,\qquad T_{2j+1}=\Theta\cdot
\]
is a lift, all higher squares reducing to $u\Theta=v\Lambda=\Theta u$. As $ker(u\cdot)=vB$, the element $\Theta$ is unique up to $vB$, so $\theta:=\mu(\Theta)\in A$ is well defined modulo $Nx^{N-1}A$. Explicitly, a telescoping computation shows that
\[
\Theta=-\sum_{l\geq 0}c_l\sum_{a=l}^{\min(N-1,\;N-2+l)}(a-l+1)\;y^az^{\,N-2+l-a}
\]
satisfies $u\Theta=v\Lambda$, the single consistency condition arising in the computation being $Nc_0=0$, once more the cocycle condition. Since $\mu$ kills every monomial of total degree $\geq N$, only the summands $l=0,1$ contribute to $\theta$, and
\[
\theta=-\binom{N}{2}\big(c_0x^{N-2}+c_1x^{N-1}\big)=-\binom{N}{2}\,x^{N-2}\lambda.
\]

\emph{Cap products.} By Theorem \ref{thm:chainmaps}, $[a]_n\cap[\lambda]_m=(-1)^{m(n-m)}\big[a\,\mu(T_{n-m}(1_B))\big]_{n-m}$, and the sign is $-1$ exactly when $m$ is odd and $n$ is even. Reading off the lifts:
\[
[a]_n\cap[\lambda]_m=
\begin{cases}
[a\lambda]_{n-m}, & m\text{ even},\\[2pt]
[a\lambda]_{n-m}, & m\text{ odd},\ n\text{ odd},\\[2pt]
-[a\theta]_{n-m}=\binom{N}{2}\big[a\,x^{N-2}\lambda\big]_{n-m}, & m\text{ odd},\ n\text{ even}.
\end{cases}
\]
All right-hand sides are independent of the choices made: changing $\lambda$ or $a$ by $Nx^{N-1}A$, or $\Theta$ by $vB$, changes them by elements of $Nx^{N-1}A$, which die in the odd-degree targets and are killed by $a\in Ann_A(Nx^{N-1})$, respectively by $\lambda\in Ann_A(Nx^{N-1})$, in the even-degree ones.

For the class $\chi:=[1]_2\in HH^2(A)$ this gives $[a]_n\cap\chi=[a]_{n-2}$: cap with $\chi$ realizes the $2$-periodicity of the resolution, an isomorphism $HH_n(A)\to HH_{n-2}(A)$ for $n\geq 3$ and the inclusion $Ann_A(Nx^{N-1})\hookrightarrow A$ for $n=2$. Cap with even cohomology classes is thus multiplication in $A$ combined with the periodicity, and the same holds for odd classes acting on odd homology. The remaining quadrant is where the example is nontrivial: over a field $k$, the pairing
\[
HH_{2i}(A)\otimes HH^{2j-1}(A)\longrightarrow HH_{2(i-j)+1}(A),\qquad [a]\otimes[\lambda]\mapsto\binom{N}{2}\big[a\,x^{N-2}\lambda\big],
\]
vanishes identically unless $k$ has characteristic $2$ and $N\equiv 2\pmod 4$, in which case it is nonzero, for example $[1]\cap[1]=[x^{N-2}]\neq 0$. Indeed, if $N$ is invertible in $k$ then $Ann_A(Nx^{N-1})=(x)$, so $x^{N-2}\lambda\in x^{N-1}A=Nx^{N-1}A$ vanishes in the target; if the characteristic of $k$ is a prime $p$ dividing $N$ then the target is $A$ and the value is $\binom{N}{2}ax^{N-2}\lambda$ with $\binom{N}{2}=\frac{N(N-1)}{2}\equiv 0$ for odd $p$, while for $p=2$ one has $\binom{N}{2}\equiv N/2\pmod 2$. For $k=\mathbb{F}_2$ and $N=2$, where $u=v$ and one may take $\Theta=\Lambda$, every cap product is plain multiplication: $[a]_n\cap[\lambda]_m=[a\lambda]_{n-m}$ for all $0\leq m\leq n$.
\end{example}

\begin{example}[Polynomial algebras and the interior product]\label{ex:poly}
Let $A=k[x_1,\dots,x_d]$ over any commutative unital ring $k$, let $V$ be the free $k$-module on $e_1,\dots,e_d$ and let
\[
P_j=A^e\otimes\Lambda^jV,\qquad d_j(1\otimes e_{i_1}\wedge\cdots\wedge e_{i_j})=\sum_{r=1}^{j}(-1)^{r-1}\,u_{i_r}\otimes e_{i_1}\wedge\cdots\widehat{e_{i_r}}\cdots\wedge e_{i_j},
\]
where $u_i=x_i\otimes 1-1\otimes x_i$. The $u_i$ form a regular sequence generating $ker(d_0)$ in $A^e=k[y_1,\dots,y_d,z_1,\dots,z_d]$: after the change of variables $w_i=y_i-z_i$ they become part of a polynomial coordinate system. Hence the Koszul complex $P_\bullet\to A$ is an $A^e$-free resolution of $A$, over every commutative $k$. The bimodule $A$ is killed by every $u_i$, hence all differentials of $A\tAe P_\bullet\cong A\otimes\Lambda^\bullet V$ and of $Hom_{A^e}(P_\bullet,A)\cong Hom_k(\Lambda^\bullet V,A)$ vanish, and
\[
HH_n(A)\cong A\otimes\Lambda^nV,\qquad HH^m(A)\cong Hom_k(\Lambda^mV,A).
\]
We write $a\otimes e_I=a\otimes e_{i_1}\wedge\cdots\wedge e_{i_n}\leftrightarrow a\,dx_{i_1}\wedge\cdots\wedge dx_{i_n}$, identifying $HH_n(A)$ with the module $\Omega^n_{A|k}$ of K\"ahler $n$-forms, and we encode a cocycle $t\in Hom_{A^e}(P_m,A)$ by the polyvector field
\[
\pi_t=\sum_{|S|=m}\lambda_S\,\partial_{s_1}\wedge\cdots\wedge\partial_{s_m},\qquad \lambda_S:=t(1\otimes e_S),\quad S=\{s_1<\cdots<s_m\},
\]
where $\partial_i=\partial/\partial x_i$; these are the Hochschild-Kostant-Rosenberg identifications. For $m=1$ the dictionary is pinned down on the bar resolution by the comparison $\Phi_1:\mathrm{Bar}(A)_1\to P_1$, $\Phi_1(1\otimes a\otimes 1)=\sum_iq_i\otimes e_i$, where $a\otimes 1-1\otimes a=\sum_iq_iu_i$ is the divided-difference expansion and $\mu(q_i)=\partial a/\partial x_i$: then $t\,\Phi_1(1\otimes a\otimes 1)=\sum_i\lambda_i\,\partial a/\partial x_i$, so $[t]\in HH^1(A)$ is the class of the derivation $X=\sum_i\lambda_i\partial_i$.

\emph{Lifts.} Fix lifts $\Lambda_S\in A^e$ with $\mu(\Lambda_S)=\lambda_S$. For disjoint subsets $K,S\subseteq\{1,\dots,d\}$ define $\varepsilon(K,S)\in\{\pm 1\}$ by $e_{K\cup S}=\varepsilon(K,S)\,e_K\wedge e_S$, each subset being wedged in increasing order. Define $T_j:P_{m+j}\to P_j$ on basis elements with $|I|=m+j$ by
\[
T_j(1\otimes e_I)=\sum_{S\subseteq I,\;|S|=m}\varepsilon(I\setminus S,\,S)\;\Lambda_S\otimes e_{I\setminus S}\,;
\]
for $m=1$ this is insertion in the last slot,
\[
T_j(1\otimes e_{i_1}\wedge\cdots\wedge e_{i_{j+1}})=\sum_{r=1}^{j+1}(-1)^{j+1-r}\,\Lambda_{i_r}\otimes e_{i_1}\wedge\cdots\widehat{e_{i_r}}\cdots\wedge e_{i_{j+1}}.
\]
Then $T_\bullet$ is a lift of $\widehat{t}$: the anchor is $d_0T_0(1\otimes e_S)=\mu(\Lambda_S)=t(1\otimes e_S)$, and for $j\geq 1$ both $d_jT_j(1\otimes e_I)$ and $T_{j-1}d_{m+j}(1\otimes e_I)$ are equal to
\[
\sum_{(S,\,p)}\eta(p,S)\;u_p\Lambda_S\otimes e_{I\setminus S\setminus\{p\}},
\]
the sum running over $S\subseteq I$ with $|S|=m$ and $p\in I\setminus S$, where $\eta(p,S)=\pm 1$ is defined by $e_I=\eta(p,S)\,e_p\wedge e_{(I\setminus p)\setminus S}\wedge e_S$: expanding $e_I=\varepsilon(I\setminus S,S)\,e_{I\setminus S}\wedge e_S$ first and then extracting $e_p$ from $e_{I\setminus S}$ computes the coefficient in $d_jT_j$, while extracting $e_p$ from $e_I$ first and then splitting off $e_S$ computes the one in $T_{j-1}d_{m+j}$; both expansions compute $\eta(p,S)$. (Insertion in the first slot, $\sum_S\varepsilon(S,I\setminus S)\,\Lambda_S\otimes e_{I\setminus S}$, satisfies instead $d_jT'_j=(-1)^mT'_{j-1}d_{m+j}$, as do the maps $s_i$ in the proof of Theorem \ref{thm:chainmaps}.)

\emph{Cap products.} Let $\omega=a\,dx_I\in\Omega^n_{A|k}\cong HH_n(A)$ and let $\pi=\pi_t$ be as above. By Theorem \ref{thm:chainmaps} and $e_{I\setminus S}\wedge e_S=(-1)^{m(n-m)}\,e_S\wedge e_{I\setminus S}$,
\[
[a\otimes e_I]\cap[t]=(-1)^{m(n-m)}\sum_{S}\varepsilon(I\setminus S,S)\,a\lambda_S\,[e_{I\setminus S}]=\sum_{S}\varepsilon(S,I\setminus S)\,a\lambda_S\,[e_{I\setminus S}].
\]
The right-hand side is the interior product: for $\pi=X_1\wedge\cdots\wedge X_m$ define
\[
(\iota_\pi\omega)(Y_1,\dots,Y_{n-m})=\omega(X_1,\dots,X_m,Y_1,\dots,Y_{n-m});
\]
since $dx_I(\partial_{j_1},\dots,\partial_{j_n})$ is the sign of the permutation taking $(j_1,\dots,j_n)$ to $I$, one has $\iota_{\partial_S}(dx_I)=\varepsilon(S,I\setminus S)\,dx_{I\setminus S}$. Therefore, under the identifications above,
\[
\omega\cap\pi=\iota_\pi\,\omega\qquad\text{for all integers }0\leq m\leq n,
\]
with no sign, while the product $\widetilde{\cap}$ of Section 4 computed through the chain map $T_\bullet$ carries the Koszul sign, $\omega\,\widetilde{\cap}\,\pi=(-1)^{m(n-m)}\iota_\pi\omega$. In particular, for a vector field $X=\sum_i\lambda_i\partial_i$,
\[
\big(a\,dx_{i_1}\wedge\cdots\wedge dx_{i_n}\big)\cap X=\sum_{r=1}^{n}(-1)^{r-1}\,a\lambda_{i_r}\,dx_{i_1}\wedge\cdots\widehat{dx_{i_r}}\cdots\wedge dx_{i_n}=\iota_X\,\omega:
\]
under the Hochschild-Kostant-Rosenberg identifications, the cap product is the contraction of differential forms by polyvector fields.
\end{example}

\begin{remark}\label{rmk:consistency}
At $(n,m)=(2,1)$ the sign in Theorem \ref{thm:chainmaps} can be tested directly against the formula of Section 1. The maps $\Psi_0=\mathrm{id}$, $\Psi_1(1\otimes e_i)=1\otimes x_i\otimes 1$ and
\[
\Psi_2(1\otimes e_i\wedge e_j)=1\otimes x_i\otimes x_j\otimes 1-1\otimes x_j\otimes x_i\otimes 1
\]
begin a comparison map $\Psi_\bullet:P_\bullet\to\mathrm{Bar}(A)_\bullet$ over the identity of $A$, since a direct check gives $d_2\Psi_2(1\otimes e_i\wedge e_j)=u_i\Psi_1(1\otimes e_j)-u_j\Psi_1(1\otimes e_i)$, so the class $[1\tAe 1\otimes e_1\wedge e_2]\in HH_2(A)$ is the bar class of
\[
c=1\tAe(1\otimes x_1\otimes x_2\otimes 1)-1\tAe(1\otimes x_2\otimes x_1\otimes 1).
\]
For the cocycle $t$ of the derivation $X=\sum_i\lambda_i\partial_i$, the formula of Section 1 gives
\[
c\cap t=\lambda_1\,[1\tAe 1\otimes x_2\otimes 1]-\lambda_2\,[1\tAe 1\otimes x_1\otimes 1]=\lambda_1[dx_2]-\lambda_2[dx_1]=\iota_X(dx_1\wedge dx_2),
\]
in agreement with Example \ref{ex:poly}, whereas $1\tAe T_1(1\otimes e_1\wedge e_2)=\lambda_2[dx_1]-\lambda_1[dx_2]=-\iota_X(dx_1\wedge dx_2)$: the sign $(-1)^{m(n-m)}$ of Theorem \ref{thm:chainmaps} is genuinely present.
\end{remark}

\medskip
\noindent\textbf{Acknowledgements}: The present work was mostly developed during a stay at the University of Leicester, UK. The author would like to express his gratitude to Claude Cibils and Sibylle Schroll for very fruitful discussions. During the preparation of this paper, the author received funds from CIMAT A.C., CONACyT and the University of Leicester.

\end{document}